\documentclass{ifacconf}

\usepackage{graphicx}      
\usepackage{natbib}        
\usepackage{amsmath,amssymb,xcolor}
\usepackage{algorithm}
\usepackage[noend]{algpseudocode}
\newcommand{\bi}{\begin{itemize}}
\newcommand{\ei}{\end{itemize}}
\newcommand{\vo}[1]{\textcolor[HTML]{000000}{\boldsymbol{#1}}}

\newcommand{\domain}[1]{\mathcal{D}_{#1}}

\usepackage{url}

\begin{document}
\begin{frontmatter}

\title{Optimal Transport Based Filtering with Nonlinear State Equality Constraints\thanksref{footnoteinfo}}

\thanks[footnoteinfo]{This research was sponsored by Air Force Office of Scientific Research, Dynamic Data Driven Applications Systems grant FA9550-15-1-0071}

\author{Niladri Das \& Raktim Bhattacharya}

\address{Department of Aerospace Engineering,
        Texas A\&M University, College Station, Texas, USA. (e-mail: niladridas,raktim@tamu.edu).}

\begin{abstract}
In this work we propose a framework to address the issue of state dependent nonlinear equality-constrained state estimation using Bayesian filtering. This framework is constructed specifically for a linear approximation of Bayesian filtering that uses the theory of Optimal Transport. 
As a part of this framework, we present three traditionally-used nonlinear equality constraint-preserving algorithms coupled with the Optimal Transport based filter: the equality-constrained Optimal Transport filter, the projected Optimal Transport filter, and the measurement-augmented Optimal Transport filter. In cases where the nonlinear equality-constraints represent an arbitrary convex manifold, we show that the re-sampling step of Optimal Transport filter, can generate initial samples for filtering, from any probability distribution function defined on this manifold. We show numerical results using our proposed framework.   
\end{abstract}


\end{frontmatter}

\section{Introduction}
The problem of Bayesian filtering for nonlinear stochastic systems, has attracted a myriad of approximation techniques owing to its intractability as the dimension of the problem increases. 
Particle filtering (PF) is one such approximation method that depends upon Markov Chain Monte-Carlo techniques and its variants, to approximate the posterior probability density function (PDF) faithfully as shown in \cite{doucet2009tutorial}. There is another class of a particle based filtering techniques. These were developed using the theory of Optimal Transport (OT). \cite{villani2008optimal} offers an excellent exposition on Optimal Transport.

The concept of Optimal Transport deals with synthesizing and analyzing transport map between two functions on some measure space, that are optimal regarding some cost function. In Bayesian filtering, it is the transport map between the prior and posterior PDF that we are interested in for some non-negative cost function. There are several approaches to design this transport map such as illustrated in \cite{reich2013nonparametric},\cite{el2012data}, \cite{taghvaei2016optimal}, and \cite{el2012bayesian}. In \cite{reich2013nonparametric}, the author proposes a linear transport map between the prior and posterior PDF, where the PDFs are expressed using samples. Here, the cost function is specified as the sum of Euclidean distances between these samples, which was extended in \cite{nil2} for cylindrical manifold, in the context of orbital dynamics. In this study, we expand the same formulation to incorporate filtering in presence of \textit{non-linear equality constraints} on the states.

A considerable group of practical systems where we require state estimation, have state constraints. For example, in dynamics and control, attitude estimation is an essential problem. Quaternion-based attitude estimation as illustrated by \cite{lefferts1982kalman} involves an equality norm constraint. State estimation of a simple pendulum with dynamics in Cartesian (non-minimal) co-ordinate system  has equality constraint on the length of the pendulum. State estimation in a one-dimensional inviscid and compressible hydrodynamic model as shown in \cite{Soares_Teixeira_2008} includes boundary conditions that are chosen such that density and energy are preserved. We can represent these conserved quantities as state constraints. Constraints arising is several other such systems can either be modeled as an inequality or inequality constraint. In this work our focus is on OT filtering for systems with nonlinear equality constraints.

In a particle based filtering framework might violate state constraints because of several reasons. Initial samples that are not generated from constrained manifold leads to prior samples that do not satisfy the constraints. Imperfect dynamics might also lead to such prior samples even if initial samples satisfy state constraints. Numerical error in the propagation also has similar effects irrespective of the choice of initial samples. Approximation of the Bayesian update also adds to this error. OT filtering also suffers from these issues when the system has inherent state constraints.

An accepted approach to deal with constraints in filters such as EKF and KF is \textit{clipping} as illustrated in \cite{simon2010kalman} where the predicted states are projected right on the boundary if outside the constraint region. The first weakness is that the constraint information is not accurately combined in the covariance update, which can contribute to a weak estimate as stated in \cite{kandepu2008applying}. For PF, we employ constraint information about the states to construct an acceptance/rejection algorithm. This is worked out by adjusting the weight calculation step, using a switch function that is 1 if the sample lie on the constraint surface and is 0 if it doesn’t. This technique allows two advantages: a) It ensures that the particles are entirely retained if they remain in the constraint region and this involves no extra computational cost. b) Because this preserves the character of Monte Carlo Sampling as we observe in \cite{hastings1970monte}, we can sample from non-Gaussian pdf. But this acceptance/rejection algorithm leads to a contraction in the sample size, which might contribute to poor estimates. With poor quality prior and non-linear constraints, all the samples might lie outside the constraint region, in which case the PF fails. An excellent survey of filtering with state constraints is presented in \cite{simon2010kalman} and \cite{yang2009kalman}.

In this study, we are interested in the state constraints represented as nonlinear equality constraints. 
In \cite{Soares_Teixeira_2008}, algorithms on nonlinear equality-constrained filtering for linear systems are extended to nonlinear systems when Unscented filtering is employed for state estimation. There are three algorithms. The first one is based on \textit{measurement-augmentation}, where the constraint equations are affixed to the existing measurement model as a perfect measurement. The second type of algorithm is based on posterior \textit{projection method} as explained in \cite{julier2007kalman}. If this posterior projection is recursively applied as feedback in the successive stages, that leads us the third algorithm that is named in \cite{Soares_Teixeira_2008} as \textit{equality-constrained} filtering.

In nonlinear filtering, Optimal Transport Filter (OTF) offers better posterior estimates than Unscented Kalman Filter (UKF) or the Ensemble Kalman Filter (EnKF) as shown in \cite{nil2}. Our main contribution is in extending the OTF where nonlinear state equality constraints are present. To the best of our knowledge, there is no prior work in OTF with state constraints. We followed the footsteps of \cite{Soares_Teixeira_2008} and \cite{teixeira2009state} and adapted the OTF framework with the three algorithms that are used for addressing nonlinear equality constraints. We further show through a numerical example, that coupling the measurement augmentation filtering with equality-constrained filtering provides the least constraint error. We also show in the context of sampling from a equality-constrained convex manifold, which is necessary to initialize OTF, the re-sampling step of OTF can generate samples for arbitrary PDF defined on this manifold. 

\section{Problem Statement}
\subsection{Bayesian Filtering}

Data assimilation algorithms that requires sequential filtering typically involves a \textit{propagation} phase  and an \textit{update}  phase. At the end of propagation, a prior ensemble is obtained by applying the state transition model to each update ensemble member produced at the previous time. In the update phase, the prior ensemble is updated to incorporating the new information available in the form of direct or indirect measurements.

We use $\boldsymbol{x}_k\in \mathbb{R}^n$ to denote the state variables and $\boldsymbol{y}_k\in \mathbb{R}^m$ is used to represent measurement variable, at time step $k$. The recursion starts from the initial distribution $p(\boldsymbol{x}_{0})$. In the prediction step, we evaluate the prior distribution of $\boldsymbol{x}_{k}$, given $\boldsymbol{y}_{1},\cdots ,\boldsymbol{y}_{k-1}$ as

\begin{align}
p(\boldsymbol{x}_{k}|\boldsymbol{y}_{1},\cdots ,&\boldsymbol{y}_{k-1})=\nonumber\\&\int{p(\boldsymbol{x}_{k}|\boldsymbol{x}_{k-1})p(\boldsymbol{x}_{k-1}|\boldsymbol{y}_{1},\cdots ,\boldsymbol{y}_{k-1})d\boldsymbol{x}_{k-1}}.\nonumber
\end{align}

where $p(\boldsymbol{x}_{k}|\boldsymbol{y}_{1},\cdots ,\boldsymbol{y}_{k-1})$ is the posterior distribution at $k^{\mathrm{th}}$ time step, $p(\boldsymbol{x}_{k}|\boldsymbol{x}_{k-1})$ encapsulates the state transition model from $k-1$ to $k$, and 
$p(\boldsymbol{x}_{k-1}|\boldsymbol{y}_{1},\cdots ,\boldsymbol{y}_{k-1})$ is the posterior state distribution at $(k-1)^{\mathrm{th}}$ time step.

In the next update state, given the new measurement $\boldsymbol{y}_{k}$ at time step $k$, the updated posterior distribution of the state $\boldsymbol{x}_{k}$ is computed using the Bayes rule as
\begin{align}
p(\boldsymbol{x}_{k}|\boldsymbol{y}_{1},\cdots,\boldsymbol{y}_{k-1},&\boldsymbol{y}_{k})\nonumber\\&=\frac{p(\boldsymbol{y}_{k}|\boldsymbol{x}_{k})p(\boldsymbol{x}_{k}|\boldsymbol{y}_{1},\cdots ,\boldsymbol{y}_{k-1})}{\int{p(\boldsymbol{y}_{k}|\boldsymbol{x}_{k})p(\boldsymbol{x}_{k}|\boldsymbol{y}_{1},\cdots ,\boldsymbol{y}_{k-1})}d\boldsymbol{x}_{k}}.
\label{bayes1}
\end{align}

\subsection{Nonlinear Equality-Constrained Bayesian Filtering}
Assume that for all $k\geq 1$, the state vectors $\boldsymbol{x}_{k}$ satisfies the equality constraint
\begin{align}
\boldsymbol{g}(\boldsymbol{x}_{k}) = \boldsymbol{d}_{k}
\label{eqc1}
\end{align}
where $\boldsymbol{g}:\mathbb{R}^n \rightarrow \mathbb{R}^s$
and $\boldsymbol{d}_k\in\mathbb{R}^s$ are known. The objective of a nonlinear equality  constrained Bayesian filtering is to evaluate eq.(\ref{bayes1}), such that it also satisfies eq.(\ref{eqc1}).

However, a closed form solution of the update equation (eq.(\ref{bayes1})) is available only for a few special cases, such as linear and Gaussian case (the Kalman filter). In ensemble or particle based filtering, where the prior and the posterior distributions are approximated by a discrete set (known as ensemble) of sample points or particles as shown in \cite{gillijns2006ensemble}, approximation techniques are available to sample from the posterior in eq.(\ref{bayes1}) . In this work we first utilize the Optimal Transport  based approximation of the Bayesian update \cite{reich2013nonparametric} and adapt it with the algorithms from \cite{teixeira2009state}, to develop a comprehensive framework for OTF with nonlinear state equality constraints.  

\section{OT Filtering}
A  typical OT problem involves two probability measures $\mu$ and $\nu$ defined on some space $X$ and a non-negative cost function $c(\cdot)$ on $X\times X$. The choices of measure spaces and cost functions are subjective and depend on the particular problem as shown in \cite{villani2008optimal}. We consider a map $\phi:\domain{\mu}\rightarrow \domain{\nu}$, where $\domain{\mu}$ denotes the support of the measure $\mu$. The objective is to transport a single element $\vo{x}$ from the support of $\mu$ to the element $\phi(\vo{x})$ in support of $\nu$, minimizing total cost which is essentially
$\inf_{\phi(\cdot)} \int c(\vo{x},\phi(\vo{x}))\mu(d\vo{x})$. Such a map $\phi^{*}$ thus obtained, is called an optimal map and  push forward of $\mu$ is $\nu$, denoted  by $\phi^{*}_{\# \mu}=\nu$.

In  \cite{reich2013nonparametric} this transport map is approximated by a linear map  $\vo{\Phi}$. The transport map is represented as a matrix  $\vo{\Phi}:=[\phi_{ij}]$  such that
 \begin{align}
 \boldsymbol{x}_{j,k}^+=\mathbb{E}(\sum_{i=1}^{N}\boldsymbol{x}_{i,k}^-\mathbb{I}_{ij})=\sum_{i=1}^{N}\boldsymbol{x}_{i,k}^-\phi_{ij}
 \label{ot1}
 \end{align}
 where $\phi_{ij}=p(\mathbb{I}_{ij})$ i.e. the probability that $\boldsymbol{x}_{i,k}^-$ is mapped to $\boldsymbol{x}_{j,k}^+$ and where the indicator $\mathbb{I}_{ij}$ takes the value 1 if $\boldsymbol{x}_{i,k}^-$ is mapped to $\boldsymbol{x}_{j,k}^+$ and 0 otherwise and $N$ denotes the sample size. 
The variable  $\boldsymbol{x}_{i,k}^-$ denotes the $i^{\mathrm{th}}$ prior sample and $\boldsymbol{x}_{i,k}^+$ denotes the $i^{\mathrm{th}}$ posterior sample, at time step $k$.  We can write eq.(\ref{ot1}) in a compact form as:
\begin{align}
\vo{X}_k^{+} &= \vo{X}_k^{-}\vo{\Phi}\nonumber\\
&=   \vo{X}_k^{-}N\vo{T}\nonumber
\end{align}
where
\textit{equally weighted} prior ensemble at time step $k$ is denoted by $\vo{X}_k^{-}\in [\vo{x}_{1,k}^{-},\vo{x}_{2,k}^{-},...,\vo{x}_{N,k}^{-}]$ and the \textit{equally weighted} posterior ensemble by $\vo{X}_k^{+}\in [\vo{x}_{1,k}^{+},\vo{x}_{2,k}^{+},...,\vo{x}_{N,k}^{+}]$.
 The matrix $\vo{T}:=[t_{ij}]$ is measure preserving, which enforces the usage of the following constraints
 \begin{equation}
 \sum_{i=1}^{N}t_{ij}=1/N,\; \sum_{j=1}^{N}t_{ij}=w_i,\text{ and }\;  t_{ij} \ge 0;
 \end{equation}
 where $w_i \propto$ likelihood-function($\vo{y}_i,\vo{x}_{i,k}^{-}$). We obtain the optimal matrix $\vo{T}$ by solving the following linear programming problem with problem size $N^2$,
\begin{align}
\vo{T}^* = &  \operatorname*{argmin}_{\vo{T}} \sum_{i=1}^N \sum_{j=1}^N t_{ij}D(\vo{x}_{i,k}^-,\hat{\vo{x}}_{j,k}^+) \\
\text{ subject to :} \nonumber \\
& \sum_{i=1}^{N}t_{ij}=1/N,  \nonumber\\
& \sum_{j=1}^{N}t_{ij}=w_i, \nonumber \\
& A_{ij} \ge 0, \nonumber
\end{align}
where $D(\vo{x}_{i,k}^-,\vo{x}_{j,k}^+)$ is the Euclidean distance between $\vo{x}_{i,k}^-$ and $\hat{\vo{x}}_{j,k}^+$. The $\hat{\vo{x}}_{j,k}^+$ denotes posterior samples with sample weights equal to $w_i$.  The equally weighted prior and unequally weighted (weighted with $w_i$) posterior have the same sample locations. Hence, the step, $\vo{X}_k^{+} = \vo{X}_k^{-}\vo{\Phi} = \vo{X}_k^{-}N\vo{T} $ can be interpreted as a \textit{re-sampling}. In a succeeding sections  we exploit this re-sampling strategy to generate samples from convex manifolds. 

\section{Constrained OT Filters}
In this section applying Optimal Transport based filtering, we show three distinct techniques for Bayesian filtering that covers non-liner equality constraints. Akin to the case of UKF for nonlinear equality constrains in \cite{Soares_Teixeira_2008}, these variations of OT filters do not ensure that the constraints will be strictly satisfied. They yield approximate results.

\subsection{Equality-Constrained OT filtering}
In this method we first generate the posterior samples $\{\vo{x}_{i,k}^{+}\}$ using normal OT filtering step. Each of these posterior samples $\{\vo{x}_{i,k}^{+}\}$ are then passed through the constraint function
\begin{align}
  \vo{\mathcal{D}}_{i,k} = \vo{g}(\vo{x}_{i,k}^{+})
\end{align}
such that $\vo{\hat{d}}_k$,$\vo{\Sigma}_{k}^{dd}$, and $\vo{\Sigma}_{k}^{xd}$ are given by
\begin{align}
  \vo{\hat{d}}_k &= \mathbb{E}[\vo{\mathcal{D}}_{i,k} ]\\
  \bar{\vo{x}}_{k}^{+} &= \mathbb{E}[\vo{x}_{i,k}^{+}]\\
  \vo{\Sigma}_{k}^{dd} &= \mathbb{E}[(\vo{\mathcal{D}}_{i,k}- \vo{\hat{d}}_k)(\vo{\mathcal{D}}_{i,k}- \vo{\hat{d}}_k)^T]\\
  \vo{\Sigma}_{k}^{xd} &= \mathbb{E}[(\vo{x}_{i,k}^{+}-\bar{\vo{x}}_{k}^{+} )(\vo{\mathcal{D}}_{i,k}- \vo{\hat{d}}_k)^T].
\end{align}
This step is different than that of UKF in \cite{teixeira2009state}  where the samples have unequal weighting. The OTF posterior samples here are all equally weighted. At this point similar to that of EnKF, we apply Kalman update on each of the samples to reduce the error with respect to the constraint equation, by calculating $\vo{K}^{p}_k$ and $\{\vo{x}_{i,k}^{p+}\}$ where the superscript $p$ denotes \textit{projection}.
They are calculated as:
\begin{align}
  \vo{K}^{p}_k &=  \vo{\Sigma}_{k}^{xd}(\vo{\Sigma}_{k}^{dd})^{-1}\\
  \vo{x}_{i,k}^{p+} &= \vo{x}_{i,k}^{+} + \vo{K}^{p}_k(\vo{\mathcal{D}}_{i,k}- \vo{\hat{d}}_k)
\end{align}
These projected samples are then used as true posterior in the subsequent propagation phase of the OT filter as a feedback. We call this filtering process \textit{nonlinear equality-constrained OT filtering} (OTNLeq).
\subsection{Projected OT filtering (OTProj)}
If we do-not use this projected samples in the subsequent time propagation of posterior samples, we end up with \textit{Projected OT filtering}.
\subsection{Measurement-Augmented OT filtering (OTMA)}
In this third technique, we utilize the concept of \textit{perfect measurements} in the context of including the equality state constraints in OT filtering. We replace the actual observation model by the augmented model
\begin{align}
  \bar{\vo{y}}_k = \bar{\vo{h}}(\vo{x}_k) + \bar{\vo{v}}_k
\end{align}
where $\bar{\vo{h}}(\vo{x}_k) = \begin{bmatrix}
\vo{h}(\vo{x}_k) \\ \vo{g}(\vo{x}_{k})\end{bmatrix}.$
The weight calculation in the OT filter now utilizes this augmented measurement model. Theoretically the pdf corresponding to the perfect measurement section should be a \textit{dirac-delta}. Practically this will inevitably lead to rejection of most of the samples by assigning zero weights to them. To avoid this issue, a tight Normal distribution is used as the corresponding measurement noise.

There is another possibility to reinforce the nonlinear equality constraints using the concepts of measurement-augmentation and projection feedback together. This is essentially OTNLeq combined with OTMA, which we name as OTNLeqMA.


\section{Initial Sample Generation: Sampling from a Constrained Manifold}
The OT filtering and its nonlinear equality constraint-addressing variants discussed in this work, requires initial samples $\{\vo{x}_{i,0}\}$ to initialize the filtering process. We might need to sample them from a PDF $\{\vo{x}_{i,0}\}\sim p_0$ with an arbitrary manifold $\mathcal{M}$ as its support. There exists different sampling techniques such as: inverse transformation method, accept-reject methods, Markov Chain Monte Carlo methods as illustrated in \cite{devroye1986sample}. We present another alternative sampling technique based on the re-sampling interpretation of the OT filtering step: $\vo{X}_k^{+} = \vo{X}_k^{-}\vo{\Phi} = \vo{X}_k^{-}N\vo{T} $. In this particular step, the equally weighted prior samples goes to equally weighted posterior samples through an intermediate step that has unequally weighted ($w_i$) posterior samples.

The idea behind using it as a sampling strategy is ensuring that intermediate posterior samples  represent our target distribution to sample from and that the unequal weights are the PDF values evaluated at each of those samples. This ensures that the final posterior samples are equal weighted random samples from our target distribution.

Suppose our target distribution is $p_t$ with a support in an arbitrary manifold $\mathcal{M}$. We assume that we have $N$ sample $\{\vo{s}_{i}^{-}\}$ from the same manifold but from a easy-to-sample PDF $p_e$. We can assume $p_e$ to be a uniform distribution $U$ such that $\vo{s}_{i}^{-}\sim U(\mathcal{M})$ . The intermediate $N$ samples $\{\hat{\vo{s}}_{i}^{+}\}$ have the exact same sample location but non-uniform weights $w_i \propto p_t(\vo{s}_{i}^{-})$. We calculate the matrix $\vo{T}$. The resulting $$\vo{s}_{i}^{+} = N\sum_{i=1}^{N}\vo{s}_{i}^{-}t_{ij}$$ are then the uniform samples drawn from $p_t$.

Although this technique opens up several avenues that deserves investigation such as :  generating smaller group of samples  in a parallel manner speeding up the process of sampling, it has a shortcoming. The transport map $\vo{\Phi} = N*\vo{T}$, renders the target samples $\{\vo{s}_{i}^{+} \}$ always as a convex combination of its initial samples $\{\vo{s}_{i}^{-}\}$. This restricts the manifold $\mathcal{M}$ to be strictly convex. A possible solution to sample from a non-convex manifold would be to approximate this manifold as a combination of smaller convex manifolds, sampling from each one of them, and combining them together.

\begin{figure}[!t]
\begin{center}
\includegraphics[width=0.3\textwidth]{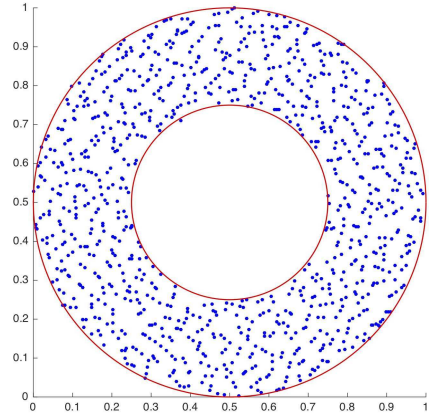}
\caption{Uniform samples generated from an annulus}
\end{center}
\label{Fig1}
\end{figure}

\begin{figure}[!t]
\begin{center}
\includegraphics[width=0.4\textwidth]{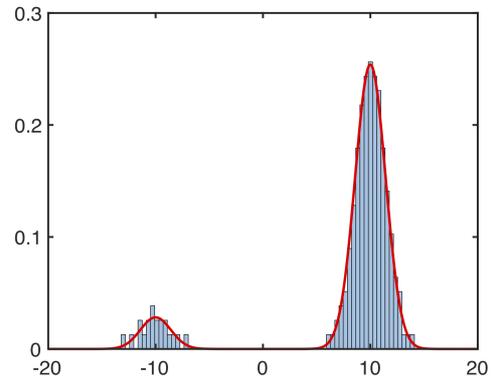}
\caption{Histogram of samples generated from Bi-modal distribution}
\end{center}
\label{Fig2}
\end{figure}

In fig.\ref{Fig1} and fig.\ref{Fig2} we show the results on implementing this sampling technique, first on an annulus with uniform distribution and then for a bimodal distribution on $\mathbb{R}$. Although the annulus is not a convex manifold, OT based sampling turns out to be an effective sampling tool. Hence, for a nonlinear system with OT filtering with nonlinear equality state constraints, we can use OT based sampling to generate the initial samples.

\section{Simulation Studies}
We consider a simple pendulum with a point mass at the end. We use a non-minimal co-ordinates $\vo{x}=[x,y]$ to denote the location of the point mass on a 2-dimensional plane. The positive x-direction point towards the right, while the positive y-direction points downwards. The kinematics of this pendulum is modeled as:
\begin{align}
\ddot{x}    &= \frac{1}{L^2}(-gxy-x(\dot{x}^2+\dot{y}^2))\\
\ddot{y}    &= \frac{1}{L^2}(gx^2-y(\dot{x}^2+\dot{y}^2))
\end{align}
where $L$ is the length of the pendulum taken to be 1 meter, and $g$ is the gravity term, which is taken to be 9.8 m/sec$^2$. Since the length of the pendulum is fixed we have $$x(t)^2+y(t)^2=L^2$$ for $t\geq 0$ as our state dependent equality constraint. We set real inital location of the pendulum point mass $\vo{x}_0=[x(0),y(0)]$ as $\vo{x}_0=[L\mathrm{cos}(30^o),L\mathrm{sin}(30^o)]$.
To perform the state estimation we assume that we can measure the location of the point mass using visual measurement with some measurement noise. The measurement model is:
\begin{align}
  \vo{y}_k = \vo{h}(\vo{x}_k) + \vo{v}_k,
\end{align}
where $\vo{h}(\vo{x}_k):=H\vo{x}_k$ with $H = \begin{bmatrix}
1 & 0 & 0 & 0 \\ 0 & 1 & 0 & 0\end{bmatrix}$. The stochastic noise variable , $\vo{v}_k \sim \mathbb{N}(\vo{0},\vo{R})$, where $\vo{R}= \begin{bmatrix}
0.01 & 0 \\  0 & 0.01 \end{bmatrix}$m$^2$. The measurements are available at every discrete time step of $\Delta t = 0.05$ secs. The constraint on the estimate at each of these time steps after assimilating the measurement is:
$$\boldsymbol{g}(\hat{\boldsymbol{x}}_k)= \hat{x}_k^2+\hat{y}_k^2,$$
where $\hat{x}_k$ and $\hat{y}_k$ denotes the position estimates of the point mass of the pendulum.

First we run OT filter without accounting for the constraints.
\begin{figure*}
\begin{center}
\includegraphics[width=\textwidth]{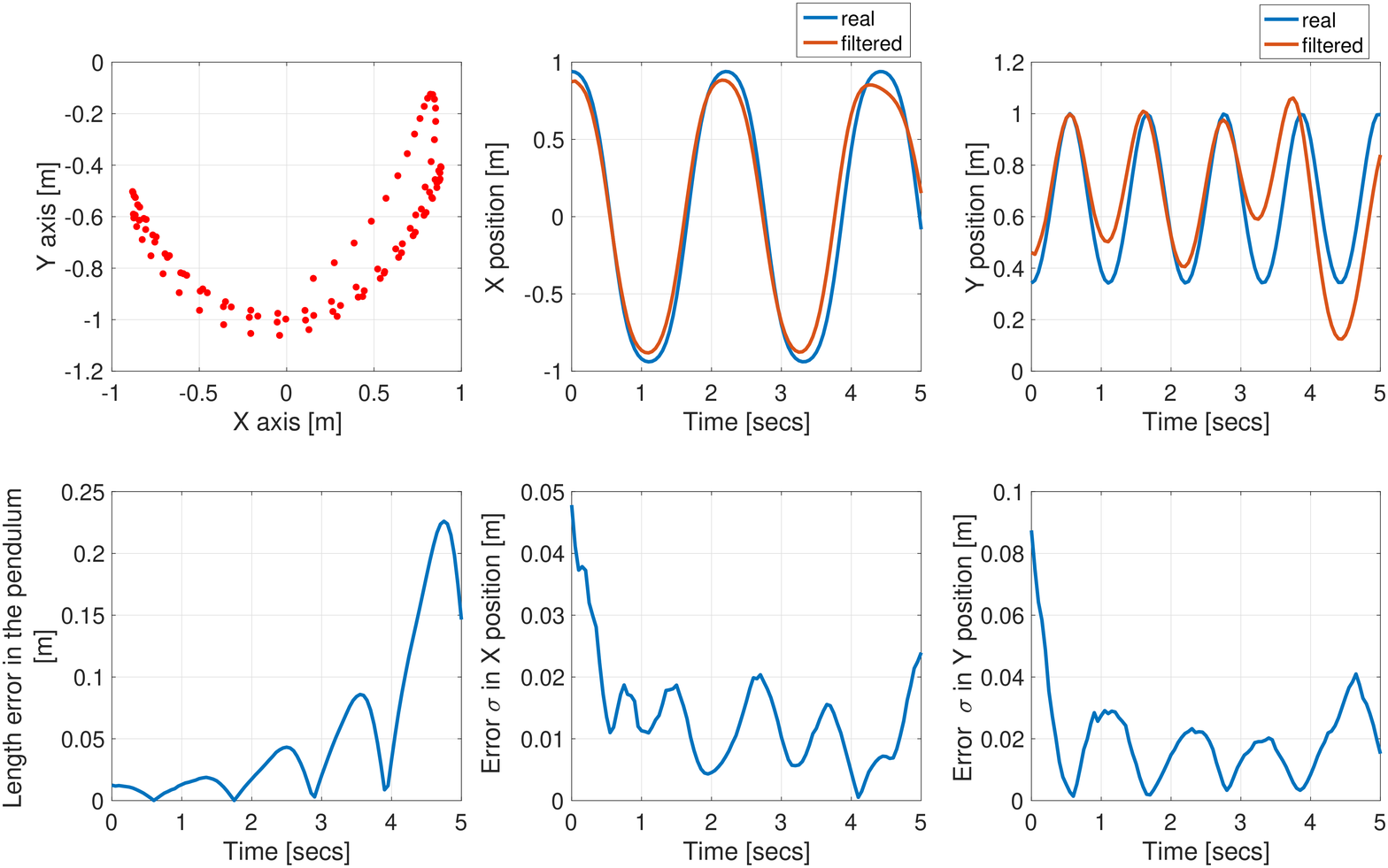}
\caption{OT filter responses without compensating for the length constraint}
\end{center}
\label{Fig3}
\end{figure*}
\begin{figure*}
\begin{center}
\includegraphics[width=\textwidth]{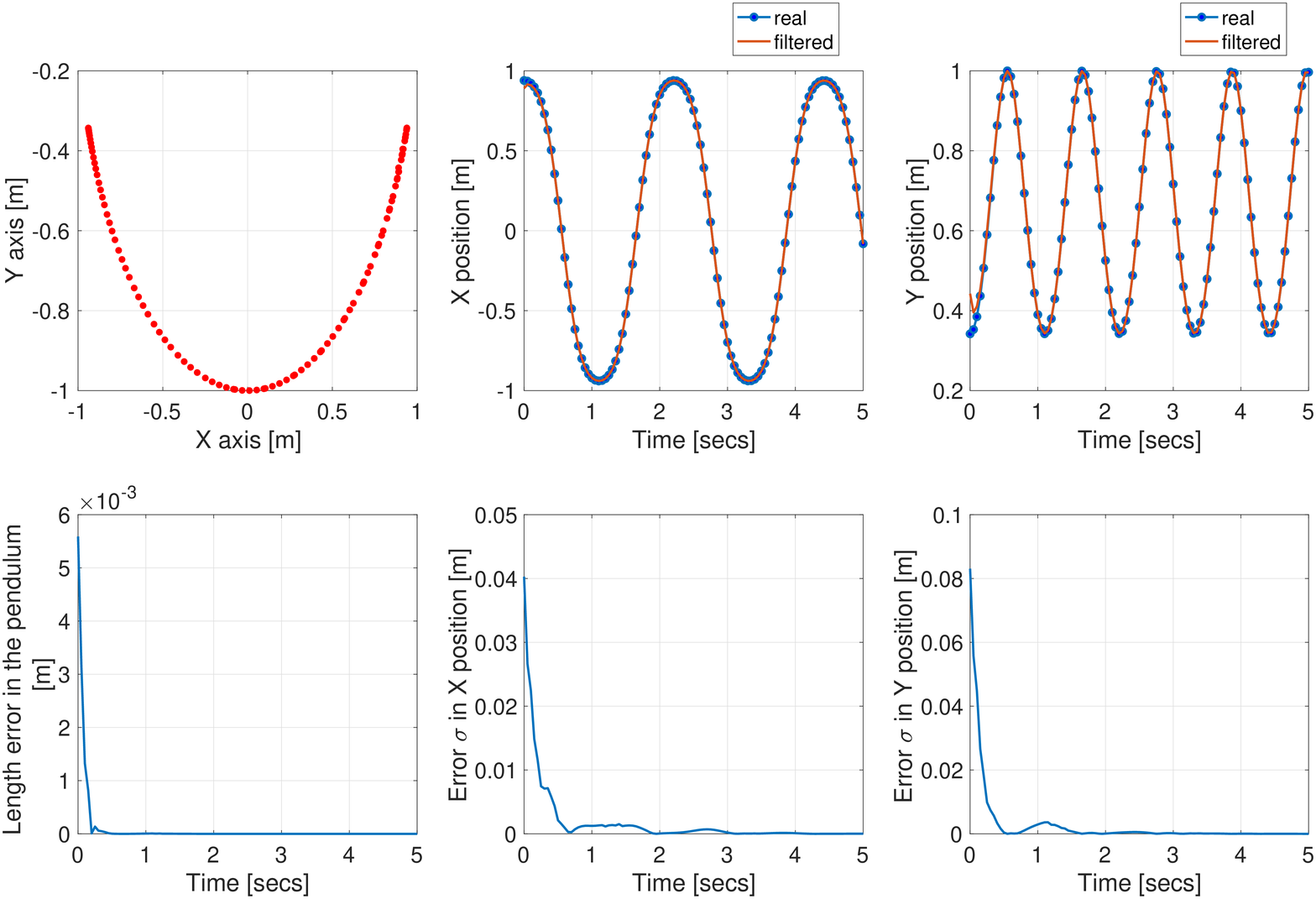}
\caption{OT filter responses after compensating for the length constraint using OTNLeqMA}
\end{center}
\label{Fig4}
\end{figure*}
In fig.\ref{Fig3} we see the response of the OT filtering when constraint preserving algorithms are not used. We notice that the state estimates starts diverging, specially the y-position estimates, from 3 seconds into the simulation. Consequently the constraint error, which is the absolute error between the true length of the pendulum and the estimated length, starts increasing with time. The trend in the state error standard deviation is also shown in this figure. 

\begin{figure}
\centering
\includegraphics[width=0.45\textwidth]{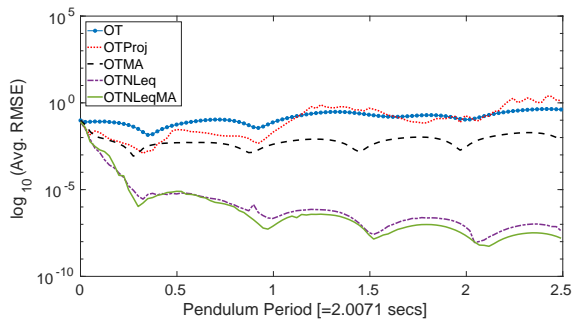}
\caption{Average RMS constraint error over 100 Monte Carlo runs}
\label{Fig5}
\end{figure}

Then, we apply all the four algorithms a). Only OTF, b). OTProj, c). OTMA, d). OTNLeq, and e). OTNLeqMA, on the same problem and run a Monte Carlo simulation with 100 runs, with different initial conditions. The average RMS error with respect to the constraint is shown in fig.\ref{Fig5}. We notice that the performance with respect to the state error constraint, of that of OTProj is comparable to using only OT filtering for state estimation. Measurement-augmentation do improve the performance but is overshadowed by that of OTNLeq and OTNLeqMA. We further notice that OTNLeqMA offers the best performance when compared to OTNLeq. This numerically shows that our new constrained filtering techniques, borne out of augmenting two existing techniques, provides state estimates with minimum error with respect to the nonlinear equality constraint. In fig.\ref{Fig4} we see the responses of the OTNLeqMA filtering algorithm.

We used MATLAB on a core i5 computer to perform our numerical experiments. The convex optimization problem needed to calculate the transport map is a linear programming problem with equality constraints, solved using \texttt{linprog()} in MATLAB.

\section{Discussions}
We addressed the nonlinear equality-constrained filtering problem for nonlinear systems when OT filtering is used for its state estimation. Three existing methodologies for equality-constrained filtering problems using KF and UKF, were coupled with OT filtering and their performances were evaluated using a numerical example. These filters were: the nonlinear equality-constrained OT filter (OTNLeq), the Projected OT filter (OTProj), and the Measurement-Augmented OT filter (OTMA). We further showed numerically, that the proposed OTNLeqMA filter provided the least constraint error compared to others. A OT based sampling technique was also studied that can sample from any PDF with an arbitrary convex manifold as its support.
\bibliography{sample}

\end{document}